# RANKS OF ELLIPTIC CURVES IN FAMILIES OF QUADRATIC TWISTS

K. RUBIN AND A. SILVERBERG

ABSTRACT. We show that the unboundedness of the ranks of the quadratic twists of an elliptic curve is equivalent to the divergence of certain infinite series.

## 1. INTRODUCTION

In this paper we reformulate the question of whether the ranks of the quadratic twists of an elliptic curve over $\mathbf{Q}$ are bounded, into the question of the whether certain infinite series converge. Our results were inspired by ideas in a paper of Gouvêa and Mazur [2].

Fix $a, b, c \in \mathbf{Z}$ such that $f(x) = x^3 + ax^2 + bx + c$ has 3 distinct complex roots, and let $E$ be the elliptic curve $y^2 = f(x)$. For $D \in \mathbf{Z} - \{0\}$, let $E^{(D)}$ be the elliptic curve $Dy^2 = f(x)$.

For every rational number $x$ which is not a root of $f(x)$, there are a unique squarefree integer $D$ and rational number $y$ such that $(x, y) \in E^{(D)}(\mathbf{Q})$. For all but finitely many $x$, the point $(x, y)$ has infinite order on the elliptic curve $E^{(D)}$. In [2], Gouvêa and Mazur count the number of $D$ that occur this way as $x$ varies, and thereby obtain lower bounds for the number of $D$ in a given range for which $E^{(D)}(\mathbf{Q})$ has positive rank.

Building on their idea, in this paper we keep track not only of the number of $D$ which occur, but also how often each $D$ occurs. The philosophy is that the greater the rank of $E^{(D)}$, the more often $D$ should occur, i.e., curves of high rank should "rise to the top". By implementing our approach, Rogers [10] found a curve of rank 6 in the family $Dy^2 = x^3 - x$.

Let

$$F(u,v) = v(u^3 + au^2v + buv^2 + cv^3) = v^4 f(u/v),$$

$$\Psi = \{(u,v) \in \mathbf{Z}^2 : \gcd(u,v) = 1 \text{ and } F(u,v) \neq 0\}.$$

We define three families of infinite series as follows.

If $n \in \mathbf{Q}^\times$, let $s(n)$ denote the square-free part of $n$, i.e., $s(n)$ is the unique square-free integer such that $n = s(n)m^2$ with $m \in \mathbf{Q}$. Note that $s(f(u/v)) = s(F(u,v))$ for all $u, v \in \mathbf{Z}$ such that $F(u,v) \neq 0$. If $\alpha$ is a non-zero rational number, and $\alpha = u/v$ with $u$ and $v$ relatively prime integers, define

$$h(\alpha) = \max\{1, \log|u|, \log|v|\}.$$

For non-negative real numbers $j$ and $k$ define the infinite sums

$$S_E(j,k) = \sum_{(u,v)\in\Psi} \frac{1}{|s(F(u,v))|^k h(u/v)^j},$$

$$R_E(j,k) = \sum_{t=1}^{\infty} \sum_{\substack{(u,v)\in\Psi \\ t^2 | F(u,v)}} \frac{t^{2k}}{|F(u,v)|^k h(u/v)^j}.$$

Further, if $d$ is a positive integer, let

$$\Omega_d = \{\alpha \in \mathbf{Z}/d^2\mathbf{Z} : f(\alpha) \equiv 0 \pmod{d^2}\}.$$

If $d$ and $d'$ are positive integers and $\alpha \in \Omega_d$, let $\omega_{\alpha,d,d'}$ be a shortest non-zero vector in the lattice

$$\mathcal{L}_{\alpha,d,d'} = \{(u,v) \in \mathbf{Z}^2 : u \equiv \alpha v \pmod{d^2} \text{ and } v \equiv 0 \pmod{(d')^2}\}.$$

(In general there will be more than one shortest vector; just choose one of them.) Define

$$Q_E(j,k) = \sum_{\substack{d,d'=1 \\ \gcd(d,d')=1}}^{\infty} \frac{(dd')^{2k}}{\max(1, \log(dd'))^j} \sum_{\substack{\alpha \in \Omega_d \\ \omega_{\alpha,d,d'} \in \Psi}} \|\omega_{\alpha,d,d'}\|^{-4k}.$$

Our main result is the following, which will be proved in §§2–4.

**Theorem 1.** *If $j$ is a positive real number, then the following are equivalent:*

(a) $\operatorname{rank}_{\mathbf{Z}} E^{(D)}(\mathbf{Q}) < 2j$ for every $D \in \mathbf{Z} - \{0\}$,
(b) $S_E(j,k)$ converges for some $k \geq 1$,
(c) $S_E(j,k)$ converges for every $k \geq 1$,
(d) $R_E(j,k)$ converges for some $k \geq 1$,
(e) $R_E(j,k)$ converges for every $k \geq 1$,
(f) $Q_E(j,k)$ converges for some $k \geq 1$,
(g) $Q_E(j,k)$ converges for every $k \geq 1$,

It follows from Theorem 1 that for many elliptic curves $E$ and for small values of $j$, $S_E(j,k)$, $R_E(j,k)$, and $Q_E(j,k)$ diverge for all real numbers $k$.

**Example 2.** Consider the case $f(x) = x^3 - x$. Here,

$$F(u,v) = uv(u+v)(u-v).$$

If $\gcd(u,v) = 1$ and $F(u,v) \neq 0$, it is not difficult to show that $s(F(u,v)) = s(u)s(v)s(u+v)s(u-v)$. In this case $S_E(j,k)$ has a particularly simple form. The family of quadratic twists $Dy^2 = x^3 - x$ has been extensively studied.

Ranks in families of twists of elliptic curves have also been studied, for example, by Heegner [5], Kramarz [7], Satgé [11], Zagier and Kramarz [15], Gouvêa and Mazur [2], Heath-Brown [3], [4], and Stewart and Top [13].

## 2. Relating $S_E(j,k)$ to twists of $E$

If $A$ is an elliptic curve over $\mathbf{Q}$, let $\hat{h}_A : A(\bar{\mathbf{Q}}) \to \mathbf{R}_{\geq 0}$ denote the canonical height function on $A(\bar{\mathbf{Q}})$. We abbreviate $\hat{h}_D = \hat{h}_{E^{(D)}}$ for squarefree integers $D$.

If $X \subset \mathbf{R}$, define
$$T_E(j,k,X) = \sum_{\substack{D \in \mathbf{Z}-0 \\ D \text{ square-free}}} |D|^{-k} \sum_{\substack{P \in E^{(D)}(\mathbf{Q}) - E^{(D)}(\mathbf{Q})_{\text{tors}} \\ x(P) \in X}} \hat{h}_D(P)^{-j}$$

where $x(P)$ is the $x$-coordinate of $P$, and define
$$S_E(j,k,X) = \sum_{(u,v) \in \Psi, u/v \in X} \frac{1}{|s(F(u,v))|^k h(u/v)^j},$$
$$R_E(j,k,X) = \sum_{t=1}^{\infty} \sum_{\substack{(u,v) \in \Psi \\ u/v \in X, t^2 | F(u,v)}} \frac{t^{2k}}{|F(u,v)|^k h(u/v)^j}.$$

Then $S_E(j,k,\mathbf{R}) = S_E(j,k)$ and $R_E(j,k,\mathbf{R}) = R_E(j,k)$ as defined in §1. Let $T_E(j,k) = T_E(j,k,\mathbf{R})$.

If $X \subset \mathbf{R}$, define

(1) $\qquad \Sigma_{D,X} = \{(u,v) \in \Psi : u/v \in X, \, v > 0, \text{ and } s(F(u,v)) = D\}.$

If $A$ is an elliptic curve over $\mathbf{Q}$, let $A_N$ denote the $N$-torsion on $A$. It is straightforward to show the following.

**Lemma 3.** *If $D$ is a square-free integer and $X \subset \mathbf{R}$, then the map*
$$\varphi_D(u,v) = \left(\frac{u}{v}, \frac{\sqrt{F(u,v)/D}}{v^2}\right)$$

*defines a bijection*
$$\varphi_D : \Sigma_{D,X} \to \{P \in E^{(D)}(\mathbf{Q}) - E_2^{(D)}(\mathbf{Q}) : x(P) \in X\}/\pm 1.$$

**Proposition 4.** *If $j, k \geq 0$ and $X \subset \mathbf{R}$, then $T_E(j,k,X)$ converges if and only if $S_E(j,k,X)$ converges.*

*Proof.* We have
$$S_E(j,k,X) = \sum_{\substack{(u,v) \in \Psi \\ u/v \in X}} |s(F(u,v))|^{-k} h(u/v)^{-j}$$
$$= 2 \sum_{D \text{ square-free}} |D|^{-k} \sum_{(u,v) \in \Sigma_{D,X}} h(u/v)^{-j}.$$

By Lemma 3,
$$T_E(j,k,X) = 2 \sum_{D \text{ square-free}} |D|^{-k} \sum_{\substack{(u,v) \in \Sigma_{D,X} \\ \varphi_D(u,v) \notin E^{(D)}(\mathbf{Q})_{\text{tors}}}} \hat{h}_D(\varphi_D(u,v))^{-j}.$$

For $(x,y) \in E^{(D)}(\mathbf{Q})$ we have
$$\hat{h}_D(x,y) = \hat{h}_E(x, \sqrt{D}y)$$

(see the hint in Exercise 8.17 on p. 239 of [12]). For $(x,y) \in E(\bar{\mathbf{Q}})$ with $x \in \mathbf{Q}$,
$$|\hat{h}_E(x,y) - \frac{1}{2}h(x)|$$

is bounded independently of $x$ and $y$ (see Theorem VIII.9.3(e) of [12]). Therefore there is a constant $C$ (independent of $u$, $v$, $D$, and $X$) such that for $(u,v) \in \Sigma_{D,X}$,

$$|\hat{h}_D(\varphi_D(u,v)) - \frac{1}{2}h(u/v)| \leq C.$$

Except for finitely many rational numbers $u/v$, we have $\frac{1}{4}h(u/v) > C$. Therefore if either $|u|$ or $|v|$ is sufficiently large, then

(2) $$\frac{1}{4}h(u/v) \leq \hat{h}_D(\varphi_D(u,v)) \leq h(u/v).$$

It follows that the convergence or divergence of $S_E(j,k,X)$ is equivalent to that of $T_E(j,k,X)$. □

If $A$ is an elliptic curve defined over $\mathbf{R}$, let $A(\mathbf{R})^0$ denote the connected component of the identity in $A(\mathbf{R})$.

**Lemma 5.** *Suppose $A$ is an elliptic curve over $\mathbf{R}$, $P_1, \ldots, P_r \in A(\mathbf{R})^0$ are $\mathbf{Z}$-linearly independent in $A(\mathbf{R})/A(\mathbf{R})_{\mathrm{tors}}$, and $U$ is an open subset of $A(\mathbf{R})^0$. Then*

$$\lim_{B \to \infty} \frac{\#\{(n_1, \ldots, n_r) \in \mathbf{Z}^r : |n_i| \leq B, \sum n_i P_i \in U\}}{(2B+1)^r} = \mu(U)$$

*where $\mu$ is a Haar measure on $A(\mathbf{R})^0$ normalized so that $\mu(A(\mathbf{R})^0) = 1$.*

*Proof.* By Satz 10 on p. 93 of [6], if $\alpha_1, \ldots, \alpha_r \in \mathbf{R}$ are $\mathbf{Z}$-linearly independent in $\mathbf{R}/\mathbf{Q}$ and $0 \leq a \leq b \leq 1$, then

$$\lim_{B \to \infty} \frac{\#\{(n_1, \ldots, n_r) \in \mathbf{Z}^r : |n_i| \leq B, a < \langle \sum n_i \alpha_i \rangle < b\}}{(2B+1)^r} = b - a$$

where $\langle z \rangle$ denotes the fractional part of a real number $z$, i. e., $0 \leq \langle z \rangle < 1$ and $z - \langle z \rangle \in \mathbf{Z}$. Since $A(\mathbf{R})^0 \cong \mathbf{R}/\mathbf{Z}$, the lemma follows easily. □

If $A$ is an elliptic curve over $\mathbf{Q}$, let

$$h_A^{\min} = \min_{\substack{P \in A(\mathbf{Q}) \\ \hat{h}_A(P) \neq 0}} \hat{h}_A(P) > 0.$$

**Proposition 6.** *Suppose $A$ is an elliptic curve over $\mathbf{Q}$, and $j$ is a positive real number. Let $r = \mathrm{rank}_{\mathbf{Z}} A(\mathbf{Q})$.*

(i) *If $r \geq 2j$ and $U$ is a nonempty open subset of the connected component of the identity in $A(\mathbf{R})$, then*

$$\sum_{P \in (A(\mathbf{Q}) - A(\mathbf{Q})_{\mathrm{tors}}) \cap U} \hat{h}_A(P)^{-j}$$

*diverges.*

(ii) *If $r < 2j$ then there exists a constant $C_j$ depending only on $j$ (and independent of $A$) such that*

$$\sum_{P \in A(\mathbf{Q}) - A(\mathbf{Q})_{\mathrm{tors}}} \hat{h}_A(P)^{-j} \leq \#A(\mathbf{Q})_{\mathrm{tors}} (h_A^{\min})^{-j} C_j.$$

*Proof.* Suppose $P_1, \ldots, P_r$ is a **Z**-basis of $A(\mathbf{Q}) \cap A(\mathbf{R})^0$ modulo torsion. The canonical height function $\hat{h}_A$ is a quadratic form on the lattice $A(\mathbf{Q})/A(\mathbf{Q})_{\text{tors}}$, and

$$\sum_{P \in A(\mathbf{Q}) - A(\mathbf{Q})_{\text{tors}}} \hat{h}_A(P)^{-j} \geq \sum_{n_1, \cdots, n_r = -\infty}^{\infty} \hat{h}(\sum n_i P_i)^{-j}.$$

By the theory of Epstein zeta functions, the latter sum diverges if $2j \leq r$. Using Lemma 5 it is now straightforward to deduce (i).

By Proposition 1(c) in IV.4.4, Vol. II of [14], there exist a positive constant $K_r$ depending only on $r$, and a **Z**-basis $P_1, \ldots, P_r$ for $A(\mathbf{Q})/A(\mathbf{Q})_{\text{tors}}$, such that for all $(n_1, \ldots, n_r) \in \mathbf{Z}^r$,

$$\hat{h}_A(\sum_{i=1}^{r} n_i P_i) \geq K_r \sum_{i=1}^{r} n_i^2 \hat{h}_A(P_i) \geq K_r h_A^{\min} \sum_{i=1}^{r} n_i^2.$$

Let $\mathcal{E}_r(j) = \sum_{0 \neq \omega \in \mathbf{Z}^r} \|\omega\|^{-2j}$. Then

$$\sum_{P \in A(\mathbf{Q}) - A(\mathbf{Q})_{\text{tors}}} \hat{h}_A(P)^{-j} \leq \#A(\mathbf{Q})_{\text{tors}} \sum_{0 \neq \omega \in \mathbf{Z}^r} (h_A^{\min})^{-j} K_r^{-j} \|\omega\|^{-2j}$$

$$= \#A(\mathbf{Q})_{\text{tors}} (h_A^{\min})^{-j} K_r^{-j} \mathcal{E}_r(j).$$

The Epstein zeta function $\mathcal{E}_r(j)$ converges if $r < 2j$ (see I.1.4, Vol. I of [14]). Thus assertion (ii) is true with $C_j = \max_{r < 2j} (K_r^{-j} \mathcal{E}_r(j))$. □

**Remark 7.** Proposition 6(ii) remains true, with the same proof, when **Q** is replaced by a number field. Proposition 6(i) remains true, with the same proof, when **Q** is replaced by a number field with a real embedding, or when **Q** is replaced by an arbitrary number field and $U$ is replaced by $A(\mathbf{C})$.

**Definition 8.** Write $e_{\max}$ (resp., $e_{\min}$) for the largest (resp., smallest) real root of $f$. We say that $X$ is *broad* if $X$ is an open subset of **R** which has nontrivial intersection with both of the intervals $(e_{\max}, \infty)$ and $(-\infty, e_{\min})$.

**Theorem 9.** *If $j$ is a positive real number, then the following are equivalent:*
  (a) $\text{rank}_{\mathbf{Z}} E^{(D)}(\mathbf{Q}) < 2j$ *for every $D \in \mathbf{Z} - \{0\}$,*
  (b) $S_E(j, k, X)$ *converges for some $k \geq 1$ and some broad $X$,*
  (c) $S_E(j, k)$ *converges for every $k \geq 1$.*

*Proof.* Fix a positive real number $j$. Clearly, (c) $\Rightarrow$ (b), by taking $X = \mathbf{R}$.

If $S_E(j, k, X)$ converges for some $k \geq 1$, and some broad $X$, then by Proposition 4, $T_E(j, k, X)$ converges as well. In particular for every square-free $D$ the inner sum

$$\sum_{\substack{P \in E^{(D)}(\mathbf{Q}) - E^{(D)}(\mathbf{Q})_{\text{tors}} \\ x(P) \in X}} \hat{h}_D(P)^{-j}$$

converges. Since $X$ is broad, the set

$$U = \{P \in E^{(D)}(\mathbf{R}) : x(P) \in X\} \cap E^{(D)}(\mathbf{R})^0$$

is nonempty. Proposition 6(i) now shows that $\text{rank}_{\mathbf{Z}} E^{(D)}(\mathbf{Q}) < 2j$. This proves that (b) $\Rightarrow$ (a).

Now suppose that $\operatorname{rank}_{\mathbf{Z}} E^{(D)}(\mathbf{Q}) < 2j$ for every $D \in \mathbf{Z} - \{0\}$. Let
$$h_D^{\min} = h_{E^{(D)}}^{\min} = \min_{\substack{P \in E^{(D)}(\mathbf{Q}) \\ \hat{h}_{E^{(D)}}(P) \neq 0}} \hat{h}_{E^{(D)}}(P).$$

By Mazur's Theorem [8], $\#E^{(D)}(\mathbf{Q})_{\text{tors}} \leq 16$. By Proposition 6(ii),
$$\sum_{P \in E^{(D)}(\mathbf{Q}) - E^{(D)}(\mathbf{Q})_{\text{tors}}} \hat{h}_D(P)^{-j} \leq 16(h_D^{\min})^{-j} C_j.$$

Therefore
$$T_E(j,k) \leq 16 C_j \sum_{\substack{D \in \mathbf{Z} - 0 \\ D \text{ square-free}}} |D|^{-k} (h_D^{\min})^{-j}.$$

It follows from Exercise 8.17c on p. 239 of [12] that there exists $D_0 > 1$, depending on $E$, such that if $|D| > D_0$, then $h_D^{\min} > \frac{1}{12}\log(|D|)$. Thus for a new constant $C'_j$,
$$T_E(j,k) \leq C'_j \left( \sum_{\substack{|D| \leq D_0 \\ D \text{ square-free}}} |D|^{-k} (h_D^{\min})^{-j} + \sum_{D > 1} |D|^{-k} (\log(|D|))^{-j} \right).$$

It follows that $T_E(j,k)$ converges if $k > 1$, or if $k = 1$ and $j > 1$. There exists a $D$ so that $\operatorname{rank}_{\mathbf{Z}} E^{(D)}(\mathbf{Q}) \geq 2$ (by [9] when the $j$-invariant of $E$ is not 0 or 1728; however, Mestre says he shows this in general in unpublished work). Therefore $j > 1$, so $T_E(j,k)$ converges. By Proposition 4, $S_E(j,k)$ converges. Therefore, (a) $\Rightarrow$ (c). □

## 3. Relating $R_E(j,k)$ and $S_E(j,k)$

**Proposition 10.** *If $k > 1/2$, $j \geq 0$, and $X \subset \mathbf{R}$, then:*
  (i) $S_E(j,k,X) \leq R_E(j,k,X) \leq \zeta(2k) S_E(j,k,X)$, *and*
  (ii) $R_E(j,k,X)$ *converges if and only if $S_E(j,k,X)$ converges.*

*Proof.* We have
$$S_E(j,k,X) = \sum_{(u,v) \in \Psi, u/v \in X} |s(F(u,v))|^{-k} h(u/v)^{-j}$$
$$\leq \sum_{t=1}^{\infty} \sum_{\substack{(u,v) \in \Psi \\ u/v \in X, t^2 | F(u,v)}} t^{2k} |F(u,v)|^{-k} h(u/v)^{-j}$$
$$= R_E(j,k,X)$$
$$\leq \sum_{n=1}^{\infty} \sum_{\substack{(u,v) \in \Psi \\ u/v \in X}} n^{-2k} |s(F(u,v))|^{-k} h(u/v)^{-j}$$
$$= \zeta(2k) S_E(j,k,X),$$

since $k > 1/2$. This is (i), and part (ii) follows immediately. □

**Corollary 11.** *If $j$ is a positive real number, then the following are equivalent:*
  (a) $\operatorname{rank}_{\mathbf{Z}} E^{(D)}(\mathbf{Q}) < 2j$ *for every $D \in \mathbf{Z} - \{0\}$,*
  (b) $R_E(j,k,X)$ *converges for some $k \geq 1$ and some broad $X$,*
  (c) $R_E(j,k)$ *converges for every $k \geq 1$.*

*Proof.* This is immediate from Proposition 10 and Theorem 9. □

## 4. Relating $Q_E(j,k)$ and $R_E(j,k)$

Let $\nu(d)$ denote the number of prime divisors of $d$. Let
$$\mathcal{S} = \{(\alpha, d, d') : d, d' \in \mathbf{Z}^+, \gcd(d, d') = 1, \alpha \in \Omega_d\}.$$

**Lemma 12.** *Suppose $(u,v) \in \Psi$, $t \in \mathbf{Z}$, and $t^2 | F(u,v)$. Then there exists a unique triple $(\alpha, d, d') \in \mathcal{S}$ such that $(u,v) \in \mathcal{L}_{\alpha,d,d'}$ and $dd' = t$.*

*Proof.* Note that $F(u,v) = v(v^3 f(u/v))$ and $v^3 f(u/v) \in \mathbf{Z}$. Since $\gcd(u,v) = 1$, we have $\gcd(v, v^3 f(u/v)) = 1$. Let
$$d = \sqrt{\gcd(t^2, v^3 f(u/v))}, \quad d' = \sqrt{\gcd(t^2, v)}, \quad \text{and} \quad \alpha = uv' \pmod{d^2},$$
where $v'$ is the inverse of $v \pmod{d^2}$. The proof is now straightforward. □

**Proposition 13.** *If $k > 1/2$ and $j \geq 0$, then $Q_E(j,k)$ converges if and only if $R_E(j,k)$ converges.*

*Proof.* It follows from Lemma 12 that
$$\{(u,v) \in \Psi : t^2 \mid F(u,v)\} = \coprod_{\substack{dd'=t \\ \gcd(d,d')=1}} \coprod_{\alpha \in \Omega_d} (\Psi \cap \mathcal{L}_{\alpha,d,d'}). \tag{3}$$

Hence if $X \subset \mathbf{R}$ we have
$$R_E(j,k,X) = \sum_{\substack{d,d'=1 \\ \gcd(d,d')=1}}^{\infty} (dd')^{2k} \sum_{\alpha \in \Omega_d} \sum_{\substack{(u,v) \in \Psi \cap \mathcal{L}_{\alpha,d,d'} \\ u/v \in X}} |F(u,v)|^{-k} h(u/v)^{-j}. \tag{4}$$

In the remainder of this proof, unless otherwise noted (by a subscript denoting additional dependence on something else), "$\ll$" and "$\gg$" mean up to a multiplicative constant that depends only on $F$, $j$, and $k$.

Suppose $(\alpha, d, d') \in \mathcal{S}$ and $\omega_{\alpha,d,d'} \in \Psi$. Then $\omega_{\alpha,d,d'}$ contributes to one of the terms in (4) when $X = \mathbf{R}$. Since $F$ has degree 4, $|F(\omega_{\alpha,d,d'})| \ll \|\omega_{\alpha,d,d'}\|^4$, so $\|\omega_{\alpha,d,d'}\|^{-4k} \ll |F(\omega_{\alpha,d,d'})|^{-k}$. Since the lattice $\mathcal{L}_{\alpha,d,d'}$ has area $(dd')^2$, Minkowski's Theorem implies that $\|\omega_{\alpha,d,d'}\| \ll dd'$, so $\log(dd')^{-j} \ll h(u/v)^{-j}$ where $\omega_{\alpha,d,d'} = (u,v)$. Therefore $Q_E(j,k) \ll R_E(j,k)$, so if $R_E(j,k)$ converges then $Q_E(j,k)$ converges.

Conversely, suppose $Q_E(j,k)$ converges. We will show that for some broad $X$, $R_E(j,k,X)$ converges. Then by Corollary 11, $R_E(j,k)$ converges as well.

Let $X$ be a broad bounded subset of $\mathbf{R}$ such that $f$ is nonzero on the closure of $X$ (for example, we could take $X = (e_{\min} - 2, e_{\min} - 1) \cup (e_{\max} + 1, e_{\max} + 2)$). Then on $X$, $|f| \gg_X 1$. Therefore if $u/v \in X$, then
$$|F(u,v)| = |v^4 f(u/v)| \gg_X |v|^4 \gg_X |u|^4,$$
the final inequality because $X$ is bounded. It follows that if $u/v \in X$ then
$$|F(u,v)| \gg_X \|(u,v)\|^4. \tag{5}$$

If $(u,v) \in \mathcal{L}_{\alpha,d,d'}$ then $(dd')^2$ divides $F(u,v)$; if further $F(u,v) \neq 0$, then
$$(dd')^2 \leq |F(u,v)| \ll \max(|u|,|v|)^4. \tag{6}$$

Thus $h(u/v) \gg \max(1, \log(dd'))$. By (4) and (5) we have $R_E(j,k,X) \ll_X R_1 + R_2$ where

$$R_1 = \sum_{\substack{d,d'=1 \\ \gcd(d,d')=1}}^{\infty} \sum_{\substack{\alpha \in \Omega_d \\ \omega_{\alpha,d,d'} \in \Psi}} \frac{(dd')^{2k}}{\max(1,\log(dd'))^j} \sum_{\substack{\omega \in \mathcal{L}_{\alpha,d,d'} \\ \omega \neq 0}} \|\omega\|^{-4k},$$

$$R_2 = \sum_{d,d'=1}^{\infty} \sum_{\substack{\alpha \in \Omega_d \\ \omega_{\alpha,d,d'} \notin \Psi}} \frac{(dd')^{2k}}{\max(1,\log(dd'))^j} \sum_{\omega \in \Psi \cap \mathcal{L}_{\alpha,d,d'}} \|\omega\|^{-4k}.$$

Exactly as in the proof of Proposition 6(ii), the theory of Epstein zeta functions shows that there is an absolute constant $C$ such that

$$\sum_{\substack{\omega \in \mathcal{L}_{\alpha,d,d'} \\ \omega \neq 0}} \|\omega\|^{-4k} \leq C \|\omega_{\alpha,d,d'}\|^{-4k}.$$

Therefore $R_1 \leq C Q_E(j,k)$, so $R_1$ converges.

It remains to show that $R_2$ converges. (Note that the terms in $R_2$ have no counterparts in $Q_E(j,k)$.) Fix positive integers $d$ and $d'$ and $\alpha \in \Omega_d$ such that $\omega_{\alpha,d,d'} \notin \Psi$. Let $t = dd'$ and let $\omega'$ be a shortest vector in $\mathcal{L}_{\alpha,d,d'} - \mathbf{Z}\omega_{\alpha,d,d'}$. Then $\{\omega_{\alpha,d,d'}, \omega'\}$ is a basis of $\mathcal{L}_{\alpha,d,d'}$,

$$\|\omega_{\alpha,d,d'}\|\|\omega'\| \gg \mathrm{Area}(\mathcal{L}_{\alpha,d,d'}) = t^2,$$

(7) $$\|\omega_{\alpha,d,d'}\| \ll \sqrt{\mathrm{Area}(\mathcal{L}_{\alpha,d,d'})} = t.$$

One can check that for every $m,n \in \mathbf{Z}$,

$$\|m\omega_{\alpha,d,d'} + n\omega'\|^2 \geq \frac{1}{2}\left(m^2\|\omega_{\alpha,d,d'}\|^2 + n^2\|\omega'\|^2\right).$$

Clearly $\Psi \cap \mathcal{L}_{\alpha,d,d'} \subset \mathcal{L}_{\alpha,d,d'} - \mathbf{Z}\omega_{\alpha,d,d'}$, so

$$\sum_{\omega \in \Psi \cap \mathcal{L}_{\alpha,d,d'}} \|\omega\|^{-4k} \leq 2 \sum_{n=1}^{\infty} \sum_{m=-\infty}^{\infty} \|m\omega_{\alpha,d,d'} + n\omega'\|^{-4k}$$

$$\ll \sum_{n=1}^{\infty} \sum_{m=0}^{\infty} (m^2 \|\omega_{\alpha,d,d'}\|^2 + n^2 t^4 \|\omega_{\alpha,d,d'}\|^{-2})^{-2k} \ll t^{-4k},$$

where the last inequality follows from (7) and a computation of the corresponding integral. Thus

$$R_2 \ll \sum_{d,d'=1}^{\infty} \sum_{\alpha \in \Omega_d} \frac{(dd')^{-2k}}{\max(1,\log(dd'))^j} \ll \sum_{d=1}^{\infty} \frac{3^{\nu(d)}}{d^{2k}} \sum_{d'=1}^{\infty} \frac{1}{d'^{2k}}$$

since $\#(\Omega_d) \ll 3^{\nu(d)}$. It is easy to see that $3^{\nu(d)} \ll_\varepsilon d^\varepsilon$ for every $\varepsilon > 0$. Therefore these sums converge, if $k > 1/2$. This completes the proof. □

**Corollary 14.** *If $j$ is a positive real number, then the following are equivalent:*
  (a) $\mathrm{rank}_{\mathbf{Z}} E^{(D)}(\mathbf{Q}) < 2j$ *for every* $D \in \mathbf{Z} - \{0\}$,
  (b) $Q_E(j,k)$ *converges for some* $k \geq 1$,
  (c) $Q_E(j,k)$ *converges for every* $k \geq 1$.

*Proof.* This is immediate from Proposition 13 and Corollary 11. □

Theorem 1 is now immediate from Theorem 9 and Corollaries 11 and 14.

## 5. Additional remarks

**Remark 15.** As in (6) and (7), each $\omega_{\alpha,d,d'}$ lies in an annulus $A_t$ of inner radius $C_1\sqrt{t}$ and outer radius $C_2 t$, with positive constants $C_1$ and $C_2$ depending only on $F$. If the lattices $\mathcal{L}_{\alpha,d,d'}$ were "random" lattices of area $t^2$ (with $F(\omega_{\alpha,d,d'}) \neq 0$) then one can compute that for large $t$, the expected value of $\frac{t^{2k}}{\|\omega_{\alpha,d,d'}\|^{4k}}$ in the annulus $A_t$ would be $\frac{1}{C_1^{4k-2} C_2^2 (2k-1) t}$. If we replace the corresponding terms of $Q_E(j,k)$ with this expected value, we obtain a "heuristic upper bound" for $Q_E(j,k)$ of

$$(8) \qquad O\left( \frac{1}{C_1^{4k}(2k-1)} \sum_{t=1}^{\infty} \frac{1}{t \log^{j-3}(t)} \right).$$

Here we have used that the number of $(\alpha, d, d') \in \mathcal{S}$ with $dd' = t$ is $O(4^{\nu(t)})$, and

$$\sum_{1 \leq t \leq x} 4^{\nu(t)} = O\left( x \log^3(x) \right).$$

The heuristic upper bound (8) correctly captures the fact that the divergence of $Q_E(j,k)$ is independent of $k$. On the other hand, the heuristic upper bound does not correctly predict the divergence of $Q_E(j,k)$. Note that (8) converges if and only if $j > 4$. However, it cannot be the case that $Q_E(j,k)$ converges for all $E$ and all $j > 4$, by Theorem 1 and the existence of elliptic curves over $\mathbf{Q}$ of rank greater than 8.

**Remark 16.** Another way of studying the "randomness" of the lattices $\mathcal{L}_{\alpha,d,d'}$ or their shortest vectors $\omega_{\alpha,d,d'}$ is as follows. For every $(\alpha, d, d') \in \mathcal{S}$, choose a random point $z_{\alpha,d,d'}$ in the annulus $A_{dd'}$. If $B, C \in \mathbf{R}^+$ define

$$\mathcal{S}_{B,C} = \{(\alpha, d, d') \in \mathcal{S} : dd' < B, \|z_{\alpha,d,d'}\| \leq C\sqrt{dd'}\}.$$

It is straightforward to compute that for fixed $C$ and large $B$,

$$(9) \qquad \text{the expected value of } \#\mathcal{S}_{B,C} \text{ is } O(\log^4(B)).$$

Now suppose that $E$ and $D$ are fixed, and $\operatorname{rank}_{\mathbf{Z}} E^{(D)}(\mathbf{Q}) = r$. Fix $r$ independent points $P_1, \ldots, P_r$ in $E^{(D)}(\mathbf{Q}) \cap E^{(D)}(\mathbf{R})^0$, and let $c = \left( \sum_i \sqrt{\hat{h}_{E^{(D)}}(P_i)} \right)^2$. As in the proof of Proposition 13, fix a broad bounded subset $X$ of $\mathbf{R}$ such that $f$ is nonzero on the closure of $X$, and for $B \in \mathbf{R}^+$ define

$$M_B = \{\sum_{i=1}^r n_i P_i : n_i \in \mathbf{Z}, |n_i| < \sqrt{\log(B)/2c}\} \cap \{P \in E^{(D)}(\mathbf{Q}) : x(P) \in X\}.$$

Suppose $P$ is a non-zero point in $M_B$. Then

$$(10) \qquad \hat{h}_D(P) \leq \log(B)/2.$$

Write $x(P) = u/v$ in lowest terms. By Lemma 3, $F(u,v) \neq 0$ and $s(F(u,v)) = D$. By Lemma 12, there is a unique triple $(\alpha, d, d') \in \mathcal{S}$ such that $(u,v) \in \mathcal{L}_{\alpha,d,d'}$ and $D(dd')^2 = F(u,v)$. Exactly as in (5), we have

$$\|\omega_{\alpha,d,d'}\| \leq \|(u,v)\| \ll_X |F(u,v)|^{1/4} = |D|^{1/4}\sqrt{dd'},$$

so

$$\|\omega_{\alpha,d,d'}\| \leq C'\sqrt{dd'} \tag{11}$$

for some constant $C'$ (depending only on $F$ and $X$). Using (6), (2), (10), and Lemma 3 we have

$$dd' = \sqrt{F(u,v)/D} \ll \max(|u|,|v|)^2 \ll B. \tag{12}$$

By Lemma 5,

$$\#M_B \gg_X \log^{r/2}(B). \tag{13}$$

It is not difficult to check that the fibers of the map from $M_B$ to $\mathcal{S}$ all have order bounded by 6 times the number of divisors of $D$, and it follows from this, (11), (12), and (13) that

$$\#\{(\alpha,d,d') \in \mathcal{S} : dd' < B, \|\omega_{\alpha,d,d'}\| \leq C'\sqrt{dd'}\} \gg_X \log^{r/2}(B). \tag{14}$$

Comparing (9) and (14) we conclude that if $\text{rank}_\mathbf{Z} E^{(D)}(\mathbf{Q}) > 8$ for at least one $D$, then the vectors $\omega_{\alpha,d,d'}$ are *not* distributed randomly in the annuli $A_{dd'}$.

**Remark 17.** The sum $Q_E(j,k)$ is very sensitive to the terms where $\omega_{\alpha,d,d'}$ lies close to the inner edge of the annulus $A_t$.

**Remark 18.** The reason for introducing $X$ in the sums is for the proof of Proposition 13 (see (5)).

**Remark 19.** By working a little harder in the proofs, one can show that Theorem 1 remains true if one replaces $Q_E(j,k)$ by a new sum where the condition $\omega_{\alpha,d,d'} \in \Psi$ in the definition of $Q_E(j,k)$ is replaced by the condition $F(\omega_{\alpha,d,d'}) \neq 0$.

**Remark 20.** Suppose we replace the cubic polynomial $f(x)$ by a polynomial of degree $d \geq 5$ (with distinct complex roots), and replace $F(u,v)$ by $v^m f(u/v)$ where $m$ is even and $m \geq d$. Then the resulting hyperelliptic curve has genus greater than one. Caporaso, Harris, and Mazur [1] conjectured that the number of rational points on curves of genus greater than one is bounded by a constant depending only on the genus of the curve. The conjecture of Caporaso-Harris-Mazur implies that the corresponding sums $S_E(j,k)$ and $R_E(j,k)$ converge for all $k > 1$ and $j \geq 0$ (since, conjecturally, $\#\Sigma_{D,\mathbf{R}}$ is bounded by a constant that is independent of $D$, where $\Sigma_{D,\mathbf{R}}$ is defined in equation (1)).

Department of Mathematics, Stanford University, Stanford, CA 94305
*E-mail address*: rubin@math.stanford.edu

Department of Mathematics, Ohio State University, Columbus, Ohio 43210
*E-mail address*: silver@math.ohio-state.edu